\documentclass[twoside]{article}
\usepackage[sc]{mathpazo}
\usepackage[T1]{fontenc}
\usepackage{microtype}
\usepackage{hyperref}
\usepackage{breakurl}
\usepackage{booktabs}
\usepackage{footnote}

\makesavenoteenv{tabular}
\usepackage{lipsum}

\usepackage[hmarginratio=1:1,top=32mm,columnsep=20pt]{geometry}
\usepackage[font=it]{caption}
\usepackage{paralist}
\usepackage{multicol}
\usepackage[toc,page]{appendix}
\usepackage{graphicx}
\usepackage{caption}
\usepackage{subcaption}
\usepackage{adjustbox}

\makeatletter
\newcommand*{\rom}[1]{\expandafter\@slowromancap\romannumeral #1@}
\makeatother
\usepackage{algorithm,algpseudocode}
\makeatletter
\def\therule{\makebox[\algorithmicindent][l]{\hspace*{.5em}\vrule height .75\baselineskip depth .25\baselineskip}}%

\newtoks\therules
\therules={}
\def\appendto#1#2{\expandafter#1\expandafter{\the#1#2}}
\def\gobblefirst#1{
  #1\expandafter\expandafter\expandafter{\expandafter\@gobble\the#1}}%
\def\LState{\State\unskip\the\therules}
\def\pushindent{\appendto\therules\therule}%
\def\popindent{\gobblefirst\therules}%
\def\printindent{\unskip\the\therules}%
\def\printandpush{\printindent\pushindent}%
\def\popandprint{\popindent\printindent}%

\algdef{SE}[WHILE]{While}{EndWhile}[1]
  {\printandpush\algorithmicwhile\ #1\ \algorithmicdo}
  {\popandprint\algorithmicend\ \algorithmicwhile}%
\algdef{SE}[FOR]{For}{EndFor}[1]
  {\printandpush\algorithmicfor\ #1\ \algorithmicdo}
  {\popandprint\algorithmicend\ \algorithmicfor}%
\algdef{S}[FOR]{ForAll}[1]
  {\printindent\algorithmicforall\ #1\ \algorithmicdo}%
\algdef{SE}[LOOP]{Loop}{EndLoop}
  {\printandpush\algorithmicloop}
  {\popandprint\algorithmicend\ \algorithmicloop}%
\algdef{SE}[REPEAT]{Repeat}{Until}
  {\printandpush\algorithmicrepeat}[1]
  {\popandprint\algorithmicuntil\ #1}%
\algdef{SE}[IF]{If}{EndIf}[1]
  {\printandpush\algorithmicif\ #1\ \algorithmicthen}
  {\popandprint\algorithmicend\ \algorithmicif}%
\algdef{C}[IF]{IF}{ElsIf}[1]
  {\popandprint\pushindent\algorithmicelse\ \algorithmicif\ #1\ \algorithmicthen}%
\algdef{Ce}[ELSE]{IF}{Else}{EndIf}
  {\popandprint\pushindent\algorithmicelse}%
\algdef{SE}[PROCEDURE]{Procedure}{EndProcedure}[2]
   {\printandpush\algorithmicprocedure\ \textproc{#1}\ifthenelse{\equal{#2}{}}{}{(#2)}}%
   {\popandprint\algorithmicend\ \algorithmicprocedure}%
\algdef{SE}[FUNCTION]{Function}{EndFunction}[2]
   {\printandpush\algorithmicfunction\ \textproc{#1}\ifthenelse{\equal{#2}{}}{}{(#2)}}%
   {\popandprint\algorithmicend\ \algorithmicfunction}%
\makeatother
\newlength\myindent
\usepackage{amsmath,amssymb}
\usepackage{blindtext}
\usepackage{varwidth}
\usepackage{lettrine}
\usepackage{abstract}

\usepackage{titlesec}
\renewcommand\thesection{\Roman{section}}
\titleformat{\section}[block]{\large\scshape\centering}{\thesection.}{1em}{}
\usepackage{fancyhdr}
\usepackage{hyperref}
\usepackage{booktabs}
\usepackage{multirow}
\usepackage{hhline}
\usepackage[english]{babel}
\usepackage{graphicx}
\usepackage{url}
\usepackage{amsmath,mathtools}

\title{\vspace{-15mm}%
	\fontsize{24pt}{10pt}\selectfont
	\textbf{A Seed-based Plant Propagation Algorithm: The Feeding Station Model}
	}	
\author{%
	\large
	\textsc{Muhammad Sulaiman\thanks{sulaiman513@yahoo.co.uk,} \thanks{Department of Mathematics, Abdul Wali Khan University Mardan KPK Pakistan.}, and Abdellah Salhi}\\[2mm]
    \normalsize	Department of Mathematical Sciences, University of Essex, Colchester CO4 3SQ, UK \\
	\vspace{-5mm}
	}
\date{}
\usepackage{amssymb}
\usepackage{amsmath,mathtools}

\begin{document}
\maketitle
\thispagestyle{fancy}
\smallskip
\noindent \textbf{Keywords:} Unconstrained global optimization, constrained optimization, engineering problems, \hspace{25mm}seed-based plant propagation

\begin{abstract}

\noindent The seasonal production of fruit and seeds resembles opening a feeding station, such as a restaurant agents/ customers will arrive at a certain rate and pick fruit (get served) at a certain rate following some appropriate processes. Therefore, dispersion follows the resource process. Modelling this process results in a search/ optimisation algorithm that used dispersion as an exploration tool that, if well captured, will find the optimum of a function over a given search space. This paper presents such an algorithm and tests it on non-trivial problems.

\end{abstract}

\section{Introduction}
A variety of plants have evolved in generous ways to propagate. Propagation through seeds is perhaps the most common of them all and one which takes advantage of all sorts of agents ranging from wind to water, to birds and animals. Beside propagation using runners, the strawberry plant uses seeds as well. These seeds are judiciously placed on the surface of a very tasty and brightly coloured fruit, the strawberries, which attract a variety of agents such as birds and animals including humans, which help the propagation.

Plants rely heavily on the dispersion of their seeds to colonise new territories and to improve their survival \cite{herrera2009plant,herrera2002seed}. There are a lot of studies and models of seed dispersion particularly for trees \cite{abrahamson1989plant,andersen1996plant,bryant1990plant,herrera2002seed,herrera2009plant}. Dispersion by wind and ballistic means are probably the most studied of all approaches \cite{glover2007understanding,yang2012flower,yang2013multi}. However, in the case of the strawberry plant, given the way the seeds stick to the surface of the fruit, Figure(1), \cite{du1996efficient}, dispersion by wind or mechanical means is very limited. Animals, however, and birds in particular are the ideal agents of dispersion \cite{krefting1949role,wenny1998directed,herrera2009plant,herrera2002seed}, in this case.

There are many biologically inspired optimization algorithms in the literature \cite{brownlee2011clever,yang2011nature}. Flower pollination algorithm (FPA) is inspired by the pollination of flowers through different agents \cite{yang2012flower}, the Swarm data clustering algorithm is inspired by pollination by bees \cite{kazemian2006swarm}, Particle Swarm Optimization (PSO) is inspired by the foraging behavior of a school of fish or a flock of birds, \cite{eberhart1995new,clerc2010particle}, Artificial Bee Colony (ABC) simulates the foraging behavior of honey bees \cite{karaboga2005idea,karaboga2008performance}, Firefly algorithm is inspired by the flashing fireflies when trying to attract a mate \cite{yang2010firefly,gandomi2011mixed}, Social Spider Optimization (SSO-C) is inspired by the cooperative behavior of social-spiders \cite{cuevas2014new}, to name a few of them.

The Plant Propagation Algorithm (PPA) also known as the strawberry algorithm was inspired by the way plants and specifically the strawberry plants propagate using runners, \cite{Salhi2010PPA,sulaiman2014engineering}. The attraction of PPA is that it can be implemented easily for all sorts of optimization problems. Moreover, it has few algorithm specific arbitrary parameters. PPA follows the principle that plants in good spots with plenty of nutrients will send many short runners. They send few long runners when in nutrient poor spots. Long runners PPA tries to explore the search space while short runners enable the algorithm to exploit the solution space well. It is necessary to make the performance of PPA better, in terms of convergence and efficiency.

In this paper we present a variant of PPA called the Seed-based Plant Propagation Algorithm the feeding station model (SbPPA). The main idea is inspired by the way frugivorous birds disperse the seeds of strawberry. The strawberry plants attract the frugivores and spread its seed for conservation in many habitats through long distances \cite{telleria2005conservation}. However, the spatial distribution of seeds depends on the availability of the strawberries on the plants and the number of visits by different agents to eat fruit.

SbPPA is tested on both unconstrained and constrained benchmark problems also used in \cite{kiran2013recombination,cuevas2014new}. Experimental results are presented in Tables 3-4 in terms of best, mean, worst and standard deviation for all algorithms. The paper is organised as follows: In Section \rom{2} we briefly introduce the feeding station model representing strawberry plants having fruits on them and the main characteristics of paths followed by different agents that disperse the seeds. Section \rom{3} presents the SbPPA in pseudo code form. The experimental settings, results and convergence graphs for different problems are given in Section \rom{4}. In  Section \rom{5} the conclusion and possible future work are given.

\section{Aspects of the Feeding Station Model of the Strawberry Plant}

Some animals and plants depend on each other to conserve their species \cite{stork1993extinction}. Thus, many plants require, for effective seed dispersal, the visits of frugivorous birds or animals according to a certain distribution, \cite{herrera2002seed,herrera2009plant,jordano2000fruits,debussche1994bird}.

Seed dispersal by different agents is also called ``seed shadow''; this shows the abundance of seeds spread locally or globally around parent plants. In this context, the strawberry feeding station model is divided in two parts: (1) The quantity of fruit or seeds available to agents, or the rate at which the agents will visit the plants, and (2) a probability density, that tells us about the service rate with which the agents are served  by the parent plants. This model tells us the quantity of seeds that is spread locally compared to that dispersed globally \cite{janzen1970herbivores,levin1976population,geritz1984efficacy,levin1984dispersal,augspurger1987wind}. There are two aspects that need to be balanced. First exploitation, which is represented by the dispersal of seeds around the parent plants. Secondly, exploration which ensures that the search space is well covered.

As a queuing system \cite{cooper1972introduction}, there are two basic components to this model: (1) the rate at which agents arrive at the strawberry plants, (2) the rate at which the agents eat fruit and leave the  plants to disperse the seeds. The agents arrive at plants in a random process. Assume that during any unit of time, whenever the fruits are available, at most one agent will arrive at a time to the plants, satisfying the orderliness condition. It is further supposed that the probability of arrivals of agents to the plants remain the same for a particular period of time. This means that the arrival rate of agents is higher when there are ripe fruit on the plants and remains the same for a further period when there is no fruit on plants; this is called stationarity condition. The arrival of one agent does not affect the rest of arrivals; this is called independence. Based on these assumptions, we conclude that the probability of arrival of $k$ agents during a cycle $t$ of fruit production by strawberry plants can be denoted by random variable $X'$, \cite{lawrence2002applied}. This can be expressed mathematically as
\begin{equation}
P(X'=k)=\frac{(\lambda t)^k e^{-\lambda t}}{k!},
\end{equation}
where $\lambda$ denotes the mean arrival rate of agents per unit time, $t$ is the length of the time interval. On the other hand, the time taken by agents in successfully eating fruit and leaving to disperse its seeds, in other words the service time for agents are expressed by a random variable, which follows the exponential probability distribution \cite{ang2004probability}. This can be expressed as follows,
  \begin{equation}
S(t)=\mu e^{-\mu t},
\end{equation}
where $\mu$ is the average number of agents that can eat fruit at time $t$. As some fruit goes to the ground around the plants after becoming fully ripe, this shows that the number of arrivals are less than the fruits available on plants. Mathematically, this can be expressed as the arrival
rate of agents is less than the fruits available on all plants, where $\lambda<\mu$.

We assume that the system is in steady state. Let $A$ denote the average number of agents on the plants, and $A_q$ the average number of agents in the queue. If we denote the average number of agents eating fruits by $\frac{\lambda}{\mu}$, then by Little's formulas \cite{little1961proof}, we have
 \begin{equation}
A=A_q+\frac{\lambda}{\mu},
\end{equation}
based on Equation (3), we need to maximize the following problem
\begin{eqnarray}
\begin{aligned}
 \hspace{-33mm}Maximize \mbox{ } A_q=A-\frac{\lambda}{\mu},
\end{aligned}
\end{eqnarray}

\hspace{10mm}subject to
\vspace{-7.5mm}
\begin{eqnarray}\hspace{-38mm}
\begin{aligned}
  &g_{1}(\lambda,\mu) = \lambda,\mu > 0,  \\
  &g_{2}(\lambda,\mu) = \lambda < \mu+1,  \\
\end{aligned}
\end{eqnarray}
where $A=10$, which represents the population size in the implementation. The simple limits on the variables are $0< \lambda,\mu\leq 100 $, After solving the problem we get $\lambda=1.1$, $\mu=0.1$ and $A_q=1$.

Moreover, frugivores may travel for a long distance to disperse seeds far away from parent SP; in doing so, they obey a L$\acute{e}$vy distribution \cite{thompson1942growth,van2007dispersal,reynolds2007free}.
\subsection{L$\acute{e}$vy distribution}

Randomization in metaheuristics is generally achieved by utilizing pseudorandom numbers, in light of some regular stochastic methodologies. L$\acute{e}$vy distributions is one of the probability density distributions for random variables. Here the random variables represent the directions of arbitrary flights by frugivores. This function of random variables ranges over real numbers with a domain called "search space".

The flight lengths of the agents served by SP, is assumed to be a heavy tailed power law distribution represented by,
\begin{equation}
L(s)\sim |s|^{-1-\beta},
\end{equation}
where $L(s)$ denotes the L$\acute{e}$vy distribution with index $\beta\in(0\mbox{,  }2)$.

L$\acute{e}$vy flights are a unique arbitrary excursions whose step lengths are drawn from (6). Another form of L$\acute{e}$vy distribution can be written as,

\begin{equation}
L(s, \gamma, \mu)=
\begin{dcases}
\sqrt{\frac{\gamma}{2\pi}}exp\left[-\frac{\gamma}{2(s-\mu)}\right] \left( \frac{1}{(s-\mu)}\right)^{\frac{3}{2}}, & 0<\mu <s< \infty \\
0 & Otherwise,
\end{dcases}
\end{equation}
this implies that
\begin{equation}
  \lim_{s\rightarrow\infty}L(s, \gamma, \mu)=\sqrt{\frac{\gamma}{2\pi}}\left( \frac{1}{s}\right)^{\frac{3}{2}},
\end{equation}
In terms of Fourier transform \cite{yang2011nature} the limiting value of $L(s)$ can be written as under,
\begin{equation}
  \lim_{s\rightarrow\infty}L(s)=\frac{\alpha\beta\Gamma(\beta)\sin(\frac{\pi\beta}{2})}{\pi|s|^{1+\beta}},
\end{equation}
where $\Gamma(\beta)$ is the Gamma function defined by
\begin{equation}
  \Gamma(\beta)=\int_0^\infty x^{\beta-1}e^{-x}dx.
\end{equation}
The steps are generated by using Mantegna's algorithm. This algorithm ensures the behaviour of L$\acute{e}$vy flights to be symmetric and stable as shown in Figure (3b).
\begin{figure}
        \centering
        \begin{subfigure}[b]{0.3\textwidth}
                \includegraphics[width=\textwidth]{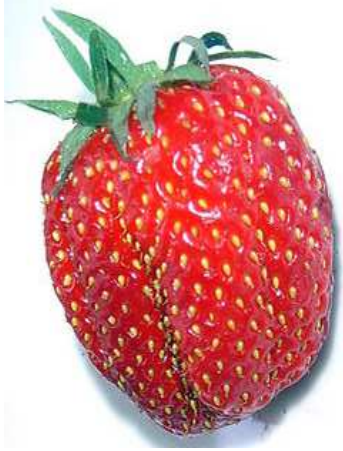}
                \caption{Strawberry fruit with seeds}
                \label{fig:gull}
        \end{subfigure}\hspace{2mm}
        ~ 
        \begin{subfigure}[b]{0.3\textwidth}
                \includegraphics[width=\textwidth]{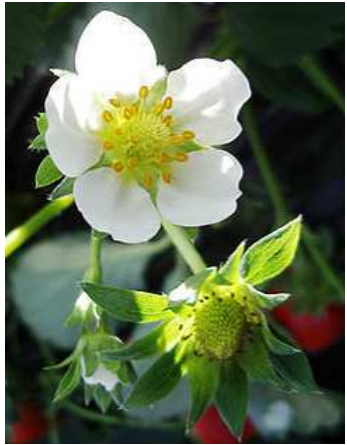}
                \caption{Strawberry garden flower\\}
                \label{fig:tiger}
        \end{subfigure}\hspace{2mm}
        ~ 
        \begin{subfigure}[b]{0.3\textwidth}
                \includegraphics[width=\textwidth]{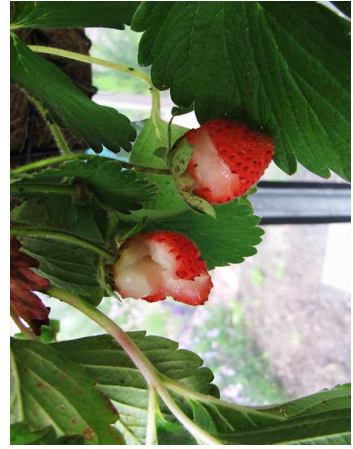}
                \caption{A fruit eaten by bird(s)}
                \label{fig:mouse}
        \end{subfigure}\\[5mm]
        \begin{subfigure}[b]{0.3\textwidth}
                \includegraphics[width=\textwidth]{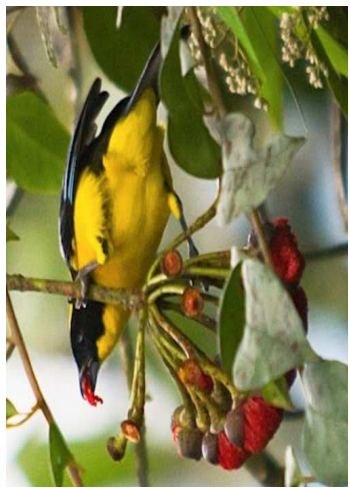}
                \caption{A bird eating strawberries}
                \label{fig:mouse}
        \end{subfigure} \hspace{9mm}
        \begin{subfigure}[b]{0.3\textwidth}
                \includegraphics[width=\textwidth]{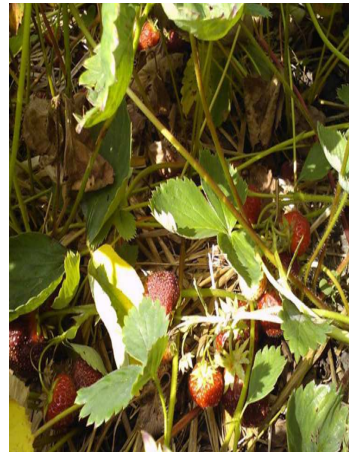}
                \caption{Strawberry plants spreading seed and sending runners around them}
                \label{fig:mouse}
        \end{subfigure}
        \caption{Strawberry plant propagation: through seed dispersion \cite{wikiStrawberry,Ruth,monacoeye,lifeisfull}}\label{fig:animals}
\end{figure}
\vspace{-3mm}
\section{Strawberry Plant Propagation Algorithm: The Feeding Station Model}
The Plant Propagation Algorithm (PPA), recently developed in \cite{Salhi2010PPA,sulaiman2014engineering}, emulates the way strawberry plants (SP) propagate by runners. Here we considered the propagation through seeds. The main objective of SbPPA is the optimal reproduction of new plants through seeds dispersion, by using different dispersal means.

We assume that the arrival of different agents to the plants for eating fruits, is according to Poisson distribution. The mean arrival rate $\lambda=1.1$, and $NP=10$ is the total number of agents in our population. Let $k=1,2,\ldots,A$ be the number of agents visiting the plants per unit time. By using these assumptions we get Figure (2) according to Equation (1).
\begin{figure}[H]
\begin{center}
\includegraphics{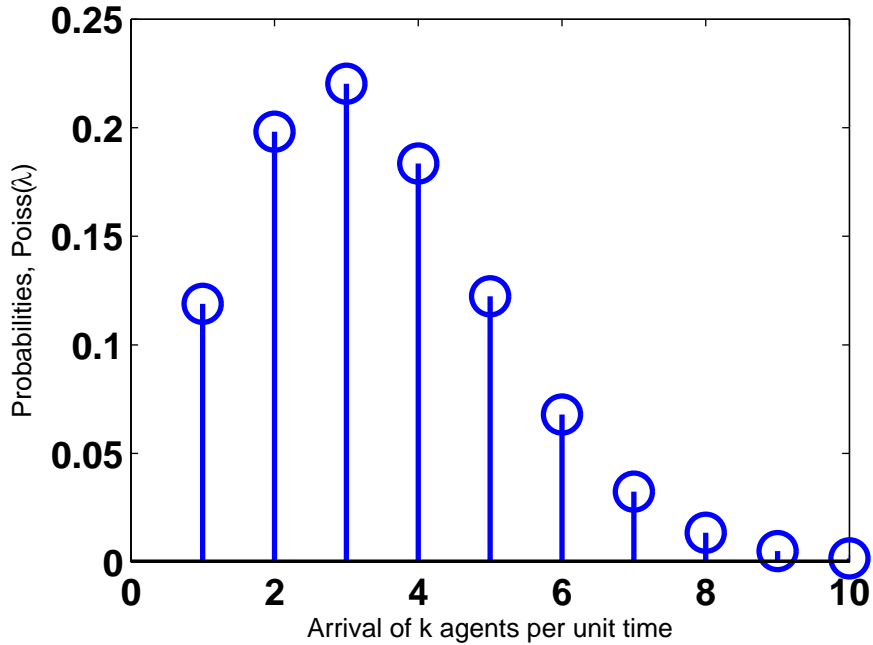}
\caption{ Agents arrival at strawberry plants to eat fruit and disperse seed}
\end{center}
\end{figure}

The probability $Poiss(\lambda)<0.05$ means that, the chances for seeds to be taken far away from SP, are lower and the propagation is supported either by runners or seeds fallen down from plants. In this case, Equation (11) below is used, which is helping the algorithm to exploit the search space,
\begin{equation}
  x^*_{i,j} = \begin{cases} x_{i,j}+\xi_j(x_{i,j}-x_{l,j}) &\mbox{if } PR\leq 0.8 \\
x_{i,j} & Otherwise, \end{cases}
\end{equation}
where $PR$ denotes the perturbation rate and it tunes the intensity of displacements by which the seeds will be dispersed locally around the SP, $x^*_{i,j}, x_{i,j}\in[a_j \mbox{ } b_j]$ are the $j^{th}$ coordinates of the seeds $X_i$ and $X^*_i$ respectively, $a_j$ and $b_j$ are the $j^{th}$ lower and upper bounds defining the search space of the problem and $\xi_j\in [-1\mbox{ } 1]$. The indices $l\mbox{ and }i$ are mutually exclusive.

On the other hand, if $Poiss(\lambda)\geq0.05$ (we choose 0.05 to give more weight to global dispersion), here the complete role of global dispersion is played by seeds, this is implemented by using the following equation,
 \begin{equation}
  x^*_{i,j} = \begin{cases} x_{i,j}+\mathrm{L_i}(x_{i,j}-\theta_j) &\mbox{if } PR\leq 0.8,\mbox{ }\theta_j\in[a_j \mbox{ }b_j]\\
x_{i,j} & Otherwise. \end{cases}
\end{equation}
Here $\mathrm{L_i}$ is a step drawn from the L$\acute{e}$vy distribution \cite{yang2011nature}, $\theta_j$ is a random coordinate within the search space. The effects on the current solutions due to perturbations applied by Equation (11) and Equation (12) are shown in Figure (3).

As mentioned in the pseudo-code of SbPPA, we first collect best solutions from the first $NP$ trial runs to form a population of potentially good solutions denoted by $pop_{best}$. The convergence rate of SbPPA, is shown in Figures (4-5), for different test problems used in our experiments. The statistical results best, worst, mean and standard deviation are calculated based on $pop_{best}$.

\begin{figure}
\includegraphics[width=\textwidth]{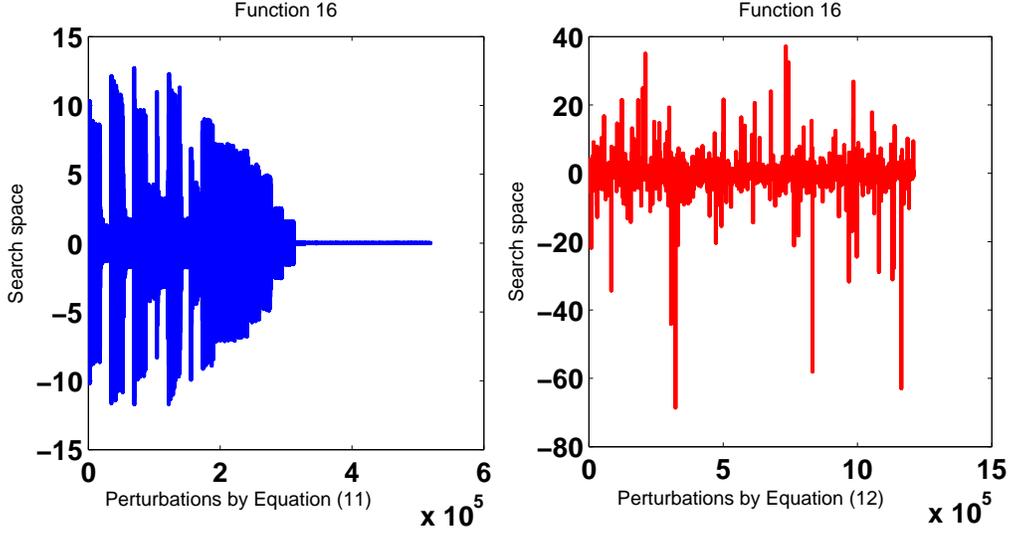}
\caption{Overall performance of SbPPA on problem 16}
\end{figure}
\newpage

The seed based propagation process of SP can be represented in the following steps:
 \begin{enumerate}
\item The dispersal of seeds or the propagation by runners in the neighbourhood of the SP, as shown in Figure ${1_e}$, is carried out either by fruit fallen from strawberry plants after they become ripe or by runners. The step lengths for this phase are calculated using Equation (11).
\item Seeds are spread globally through frugivores, as shown in Figure $1_{c,d}$. The step lengths for those travelling agents are drawn from the L$\acute{e}$vy distribution.
\item The probabilities, $Poiss(\lambda)$, that a certain amount $k$ of agents will arrive to SP to eat fruits and disperse it, is used as a switch between global and local search.
\end{enumerate}

For implementation, we assume that each SP produces one fruit, and each fruit is assumed to have one seed, we mean by a solution $X_i$ the position of the $i^{th}$ seed to be dispersed. The number of seeds in the population is denoted by $NP$. Initially we generate a random population of $NP$ seeds using Equation (13),
\begin{equation}x_{i,j}=a_j+(b_j-a_j)\eta_j, j=1,...,n\end{equation}
\noindent where $x_{i,j}\in[a_j,b_j]$ is the $j^{th}$ entry of solution $X_i$, $a_j$ and $b_j$ are the $j^{th}$ coordinates of the bounds describing the search space of the problem and $\eta_j\in (0,1)$. This means $X_{i}=[x_{i,j}], \mbox{ for } j=1,...,n$ represents the position of the $i^{th}$ seed in population $pop$.

\begin{algorithm}\label{MPPA}
\caption{\textbf{Seed-based Plant Propagation Algorithm (SbPPA): The Feeding Station Model}}
\begin{algorithmic}[1]
\LState Initialize: $g_{max}\leftarrow$ maximum number of generations, $max_{eval}\leftarrow$ maximum function evaluations, $r\leftarrow$ counter for trial runs
\LState Set $r=1$

\If {$r\leq NP$}
\LState \begin{varwidth}[t]{\linewidth}Create a random population of seeds $pop=\{X_i\mid i=1,2,...,NP\}$, using Equation (13) \par
                                       and add the best solutions from each trial run, in $pop_{best}$.
         \end{varwidth}
 \LState Evaluate the population.
\EndIf
\While {$r>NP$}
\LState Use population $pop_{best}$.
\EndWhile
\LState Set $ngen=1$,
\While {($ngen$ $<$ $g_{max}$) $\textbf{or}$ ($n_{eval} < max_{eval}$)}
\For {$i=1$ to $NP$}
 \If {$Poiss(\lambda)_i \geq 0.05 $}, \Comment (Global or local seed dispersion)
 \For {$j=1$ to $n$}   \Comment ($n$ is number of dimensions)
 \If {rand$ \leq PR $}, \Comment (PR=Perturbation rate)
 \LState Update the current entry according to Equation (12)
 \EndIf
 \EndFor
 \Else
 \For {$j=1$ to $n$}
 \If {rand$ \leq PR $},
 \LState Update the current entry according to Equation (11)
 \EndIf
 \EndFor
\EndIf
\EndFor
\EndWhile
\LState\emph{\textbf{Return}:} Update current population.
\end{algorithmic}
\end{algorithm}

\section{Experimental Setting And Discussion}
In our experiments we test SbPPA against other state-of-the-art algorithms. Our set of test problems include benchmark constrained and unconstrained optimization problems \cite{suganthan2005problem,liang2006problem,cuevas2014new}. The results are compared in terms of best, worst, mean and standard deviations obtained by SbPPA, ABC \cite{karaboga2005idea,karaboga2011modified}, PSO \cite{he2007hybrid}, FF \cite{gandomi2011mixed}, HPA \cite{kiran2013recombination} and SSO-C \cite{cuevas2014new}. The detailed descriptions of these problems are given in Appendix \rom{1}. The significance of results are shown according to the following notations:
\begin{itemize}
  \item (+) when SbPPA is better
  \item (-) when SbPPA is worse
  \item ($\approx$) when the results are approximately same as SbPPA.
\end{itemize}
\subsection{Parameter Settings}
The parameter settings are give in Table 1-2:
{\fontsize{12}{12}
\selectfont
\begin{table}[htp]
  \centering
  \footnotesize\setlength{\tabcolsep}{2.5pt}
  \begin{minipage}{15.4cm}
  \caption{Parameters used for each algorithm for solving unconstrained global optimization problems $f_{1}-f_{10}$, All experiments are repeated 30 times.}
   \begin{tabular}{l@{\hspace{6pt}} *{4}{l}}
    \toprule
    PSO \cite{eberhart1995new,kiran2013recombination} & ABC \cite{karaboga2005idea,kiran2013recombination}   & HPA \cite{kiran2013recombination}   & SbPPA \\
    \midrule
    M=100 & SN=100 & Agents=100  & NP=10 \\
    $G_{max}=\frac{(Dimension\times20,000)}{NP}$ & MCN=$\frac{(Dimension\times20,000)}{NP}$ & Iteration number=$\frac{(Dimension\times20,000)}{NP}$  & Iteration number=$\frac{(Dimension\times20,000)}{NP}$ \\[2mm]
    $c_1=2$ & MR=0.8 & $c_1=2$ & PR=0.8, $Poiss(\lambda)=0.05$ \\[2mm]
    $c_2=2$ &    limit=$\frac{(SN\times dimension)}{2}$  & $c_2=2$&$k=1,2,\ldots,A$  \\[2mm]
                                        W= $\frac{(G_{max}-iteration_{index})}{G_{max}}$  &-& limit=$\frac{(SN\times dimension)}{2}$ &  $\lambda=1.1$    &\\[2mm]
    -&  -     &  W= $\frac{(G_{max}-iteration_{index})}{G_{max}}$ &- \\
    \bottomrule
    \end{tabular}%
    \end{minipage}
  \label{tab:addlabel}%
  \end{table}}

{\fontsize{12}{12}
\selectfont
\begin{table}[htp]
  \centering
  \footnotesize\setlength{\tabcolsep}{2.5pt}
  \begin{minipage}{15cm}
  \caption{Parameters used for each algorithm for solving constrained optimization problems $f_{11}-f_{18}$, \\ All experiments are repeated 30 times.}
   \begin{tabular}{l@{\hspace{6pt}} *{5}{l}}
    \toprule
    PSO \cite{he2007hybrid}  & ABC \cite{karaboga2011modified}   & FF \cite{gandomi2011mixed}   & SSO-C \cite{cuevas2014new}& SbPPA \\
    \midrule
    M=250 & sn=40 & Fireflies=25 & N=50  & NP=10 \\[2mm]
    $G_{max}=300$ & MCN=6000 & Iteration number= 2000 & Iteration number=500 & Iteration number=2400 \\[2mm]
    $c_1=2$ & MR=0.8 & q=1.5 & PF=0.7 & PR=0.8, $Poiss(\lambda)=0.05$  \\[2mm]
    $c_2=2$ &     -  & $\alpha=0.001$ &   -    &$k=1,2,\ldots,A$  \\[2mm]
    Weight factors= 0.9 to 0.4 &  -     &  -     &   -    & $\lambda=1.1$ \\
    \bottomrule
    \end{tabular}%
    \end{minipage}
  \label{tab:addlabel}%
  \end{table}
}

{\fontsize{12}{12}
\selectfont
\begin{table}[htp]
  \centering
  \footnotesize\setlength{\tabcolsep}{2.5pt}
  \begin{minipage}{9.5cm}
  \caption{Results obtained by SbPPA, HPA, PSO and ABC. All problems in this table are unconstrained.}
   \begin{tabular}{l@{\hspace{6pt}} *{7}{l}}
     \toprule
    Fun   & Dim   & Algorithm & Best  & Worst & Mean  & SD \\
    \midrule
    1     & 4     & ABC   & (+)\mbox{ }0.0129 & (+)\mbox{ }0.6106 &(+)\mbox{ } 0.1157 &(+) \mbox{ }0.111 \\
          &       & PSO   & (-)\mbox{ }6.8991E-08 & (+)\mbox{ }0.0045 & (+)\mbox{ }0.001 & (+)\mbox{ }0.0013 \\
          &       & HPA   & (+)\mbox{ }2.0323E-06 & (+)\mbox{ }0.0456 & (+)\mbox{ }0.009 &(+) \mbox{ }0.0122 \\
          &       & SbPPA & 1.08E-07 & 7.05E-06 & 3.05E-06 & 3.14E-06 \\[0.7mm]
    2     & 2     & ABC   &(+)\mbox{ }1.2452E-08 &(+)\mbox{ }8.4415E-06 & (+)\mbox{ }1.8978E-06 & (+)\mbox{ }1.8537E-06 \\
          &       & PSO   &($\approx$)\mbox{ }0     & ($\approx$)\mbox{ }0     & ($\approx$)\mbox{ }0     & ($\approx$)\mbox{ }0 \\
          &       & HPA   &($\approx$)\mbox{ }0     &($\approx$)\mbox{ }0     &($\approx$)\mbox{ }0     &($\approx$)\mbox{ }0 \\
          &       & SbPPA & 0     & 0     & 0     & 0 \\[0.7mm]
    3     & 2     & ABC   & ($\approx$)\mbox{ }0     & (+)\mbox{ }4.8555E-06 & (+)\mbox{ }4.1307E-07 & (+)\mbox{ }1.2260E-06 \\
          &       & PSO   & ($\approx$)\mbox{ }0     & (+)\mbox{ }3.5733E-07 & (+)\mbox{ }1.1911E-08 & (+)\mbox{ }6.4142E-08 \\
          &       & HPA   & ($\approx$)\mbox{ }0     & ($\approx$)\mbox{ }0     & ($\approx$)\mbox{ }0     & ($\approx$)\mbox{ }0 \\
          &       & SbPPA & 0     & 0     & 0     & 0 \\[0.7mm]
    4     & 2     & ABC   & ($\approx$)\mbox{ }-1.03163 & ($\approx$)\mbox{ }-1.03163 & ($\approx$)\mbox{ }-1.03163 & ($\approx$)\mbox{ }0 \\
          &       & PSO   & ($\approx$)\mbox{ }-1.03163 & ($\approx$)\mbox{ }-1.03163 & ($\approx$)\mbox{ }-1.03163 & ($\approx$)\mbox{ }0 \\
          &       & HPA   & ($\approx$)\mbox{ }-1.03163 & ($\approx$)\mbox{ }-1.03163 & ($\approx$)\mbox{ }-1.03163 & ($\approx$)\mbox{ }0 \\
          &       & SbPPA & -1.031628 & -1.031628 & -1.031628 & 0 \\[0.7mm]
    5     & 6     & ABC   & ($\approx$)\mbox{ }-50.0000 & ($\approx$)\mbox{ }-50.0000 & ($\approx$)\mbox{ }-50.0000 & (-)\mbox{ }0 \\
          &       & PSO   & ($\approx$)\mbox{ }-50.0000 & ($\approx$)\mbox{ }-50.0000 & ($\approx$)\mbox{ }-50.0000 & (-)\mbox{ }0 \\
          &       & HPA   & ($\approx$)\mbox{ }-50.0000 & ($\approx$)\mbox{ }-50.0000 & ($\approx$)\mbox{ }-50.0000 & (-)\mbox{ }0 \\
          &       & SbPPA & -50.0000 & -50.0000 & -50.0000 & 5.88E-09 \\[0.7mm]
    6     & 10    & ABC   & (+)\mbox{ }-209.9929 & (+)\mbox{ }-209.8437 & (+)\mbox{ }-209.9471 & (+)\mbox{ }0.044 \\
          &       & PSO   & ($\approx$)\mbox{ }-210.0000 & ($\approx$)\mbox{ }-210.0000 & ($\approx$)\mbox{ }-210.0000 &(-)\mbox{ }0 \\
          &       & HPA   & ($\approx$)\mbox{ }-210.0000 & ($\approx$)\mbox{ }-210.0000 & ($\approx$)\mbox{ }-210.0000 & (+)\mbox{ }1 \\
          &       & SbPPA & -210.0000 & -210.0000 & -210.0000 & 4.86E-06 \\[0.7mm]
    7     & 30    & ABC   & (+)\mbox{ }2.6055E-16 & (+)\mbox{ }5.5392E-16 & (+)\mbox{ }4.7403E-16 & (+)\mbox{ }9.2969E-17 \\
          &       & PSO   & ($\approx$)\mbox{ }0     & ($\approx$)\mbox{ }0     & ($\approx$)\mbox{ }0     & ($\approx$)\mbox{ }0 \\
          &       & HPA   & ($\approx$)\mbox{ }0     & ($\approx$)\mbox{ }0     &($\approx$)\mbox{ }0     & ($\approx$)\mbox{ }0 \\
          &       & SbPPA & 0     & 0     & 0     & 0 \\[0.7mm]
    8     & 30    & ABC   & (+)\mbox{ }2.9407E-16 & (+)\mbox{ }5.5463E-16 & (+)\mbox{ }4.8909E-16 & (+)\mbox{ }9.0442E-17 \\
          &       & PSO   & ($\approx$)\mbox{ }0     & ($\approx$)\mbox{ }0     & ($\approx$)\mbox{ }0     & ($\approx$)\mbox{ }0 \\
          &       & HPA   & ($\approx$)\mbox{ }0     & ($\approx$)\mbox{ }0     & ($\approx$)\mbox{ }0     & ($\approx$)\mbox{ }0 \\
          &       & SbPPA & 0     & 0     & 0     & 0 \\[0.7mm]
    9    & 30    & ABC   & ($\approx$)\mbox{ }0     & (+)\mbox{ }1.1102E-16 & (+)\mbox{ }9.2519E-17 & (+)\mbox{ }4.1376E-17 \\
          &       & PSO   & ($\approx$)\mbox{ }0     & (+)\mbox{ }1.1765E-01 & (+)\mbox{ }2.0633E-02 & (+)\mbox{ }2.3206E-02 \\
          &       & HPA   & ($\approx$)\mbox{ }0     & ($\approx$)\mbox{ }0     & ($\approx$)\mbox{ }0     & ($\approx$)\mbox{ }0 \\
          &       & SbPPA & 0     & 0     & 0     & 0 \\[0.7mm]
    10    & 30    & ABC   & (+)\mbox{ }2.9310E-14 & (+)\mbox{ }3.9968E-14 & (+)\mbox{ }3.2744E-14 & (+)\mbox{ }2.5094E-15 \\
          &       & PSO   & ($\approx$)\mbox{ }7.9936E-15 & (+)\mbox{ }1.5099E-14 & (-)\mbox{ }8.5857E-15 & (+)\mbox{ }1.8536E-15 \\
          &       & HPA   & ($\approx$)\mbox{ }7.9936E-15 & (+)\mbox{ }1.5099E-14 & (+)\mbox{ }1.1309E-14 & (+)\mbox{ }3.54E-15 \\
          &       & SbPPA & 7.994E-15 & 7.99361E-15 & 7.994E-15 & 7.99361E-15 \\
    \bottomrule
  \end{tabular}%
    \end{minipage}
  \label{tab:CF01}%
  \end{table}
}

\begin{figure}
\hspace{-21mm}
\includegraphics[width=19cm]{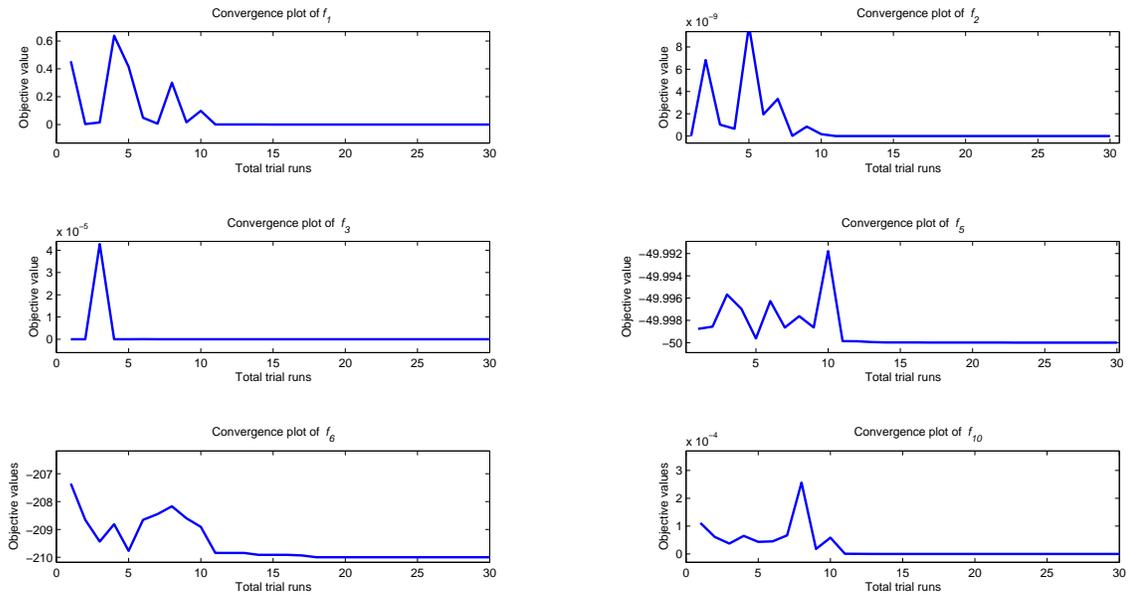}
\caption{Performance of SbPPA on unconstrained global optimization problems}

\end{figure}
{\fontsize{12}{12}
\selectfont
\begin{table}[htp]
  \centering
  \footnotesize\setlength{\tabcolsep}{2.5pt}
  \begin{minipage}{13cm}
  \caption{Results obtained by SbPPA, PSO, ABC, FF and SSO-C. All problems in this table are standard constrained optimization problems}
   \begin{tabular}{l@{\hspace{10pt}} *{8}{l}}
    \toprule
    Fun   & Fun Name & Optimal & Algorithm & Best  & Mean  & Worst & SD \\
    \midrule
    11  & CP1   & -15   & PSO   & ($\approx$)\mbox{ }-15   & ($\approx$)\mbox{ }-15   & ($\approx$)\mbox{ }-15   & (-)\mbox{ }0 \\
          &       &       & ABC   & ($\approx$)\mbox{ }-15   &($\approx$)\mbox{ }-15   & ($\approx$)\mbox{ }-15   & (-)\mbox{ }0 \\
          &       &       & FF    & (+)\mbox{ }14.999 & (+)\mbox{ }14.988 & (+)\mbox{ }14.798 & (+)\mbox{ }6.40E-07 \\
          &       &       & SSO-C & ($\approx$)\mbox{ }-15   & ($\approx$)\mbox{ }-15   & ($\approx$)\mbox{ }-15   & (-)\mbox{ }0 \\
          &       &       & SbPPA & -15   & -15   & -15   & 1.95E-15 \\[0.7mm]
   12  & CP2   & -30665.539 & PSO   & ($\approx$)\mbox{ }-30665.5 & (+)\mbox{ }-30662.8 & (+)\mbox{ }-30650.4 & (+)\mbox{ }5.20E-02 \\
          &       &       & ABC   & ($\approx$)\mbox{ }-30665.5 & (+)\mbox{ }-30664.9 &(+)\mbox{ }-30659.1 & (+)\mbox{ }8.20E-02 \\
          &       &       & FF    & ($\approx$)\mbox{ }-3.07E+04 & (+)\mbox{ }-30662 & (+)\mbox{ }-30649 & (+)\mbox{ }5.20E-02 \\
          &       &       & SSO-C & ($\approx$)\mbox{ }-3.07E+04 & ($\approx$)\mbox{ }-30665.5 & (+)\mbox{ }-30665.1 & (+)\mbox{ }1.10E-04 \\
          &       &       & SbPPA & -30665.5 & -30665.5 & -30665.5 & 2.21E-06 \\[0.7mm]
  13  & CP3   & -6961.814 & PSO   & (+)\mbox{ }-6.96E+03 & (+)\mbox{ }-6958.37 & (+)\mbox{ }-6942.09 & (+)\mbox{ }6.70E-02 \\
          &       &       & ABC   & (-)\mbox{ }-6961.81 & (+)\mbox{ }-6958.02 & (+)\mbox{ }-6955.34 & (-)\mbox{ }2.10E-02 \\
          &       &       & FF    & (+)\mbox{ }-6959.99 & (+)\mbox{ }-6.95E+03 & (+)\mbox{ }-6947.63& (-)\mbox{ }3.80E-02 \\
          &       &       & SSO-C & (-)\mbox{ }-6961.81 & (+)\mbox{ }-6961.01 & (+)\mbox{ }-6960.92 & (-)\mbox{ }1.10E-03 \\
          &       &       & SbPPA & -6961.5 & -6961.38 & -6961.45 & 0.043637 \\[0.7mm]
   14   & CP4   & 24.306 & PSO   & (-)\mbox{ }24.327 & (+)\mbox{ }2.45E+01 & (+)\mbox{ }24.843 & (+)\mbox{ }1.32E-01 \\
          &       &       & ABC   & (+)\mbox{ }24.48 & (+)\mbox{ }2.66E+01 & (+)\mbox{ }28.4  & (+)\mbox{ }1.14 \\
          &       &       & FF    & (-)\mbox{ }23.97 & (+)\mbox{ }28.54 & (+)\mbox{ }30.14 & (+)\mbox{ }2.25 \\
          &       &       & SSO-C & (-)\mbox{ }24.306 & (-)\mbox{ }24.306 & (-)\mbox{ }24.306 & (-)\mbox{ }4.95E-05 \\
          &       &       & SbPPA & 24.34442 & 24.37536 & 24.37021 & 0.012632 \\[0.7mm]
   15  & CP5   & -0.7499 & PSO   & ($\approx$)\mbox{ }-0.7499 & (+)\mbox{ }-0.749 & (+)\mbox{ }-0.7486 & (+)\mbox{ }1.20E-03 \\
          &       &       & ABC   & ($\approx$)\mbox{ }-0.7499 & (+)\mbox{ }-0.7495 & (+)\mbox{ }-0.749 & (+)\mbox{ }1.67E-03 \\
          &       &       & FF    & (+)\mbox{ }-0.7497 & (+)\mbox{ }-0.7491 & (+)\mbox{ }-0.7479 & (+)\mbox{ }1.50E-03 \\
          &       &       & SSO-C & ($\approx$)\mbox{ }-0.7499 & ($\approx$)\mbox{ }-0.7499 & ($\approx$)\mbox{ }-0.7499 & (-)\mbox{ }4.10E-09 \\
          &       &       & SbPPA & 0.7499 & 0.749901 & 0.7499 & 1.66E-07 \\[0.7mm]
   16   & Spring & Not Known & PSO   & (+)\mbox{ }0.012858 & (+)\mbox{ }0.014863 & (+)\mbox{ }0.019145 & (+)\mbox{ }0.001262 \\
          &Design       &       & ABC   & ($\approx$)\mbox{ }0.012665 & (+)\mbox{ }0.012851 & (+)\mbox{ }0.01321 & (+)\mbox{ }0.000118 \\
          &Problem       &       & FF    & ($\approx$)\mbox{ }0.012665 & (+)\mbox{ }0.012931 & (+)\mbox{ }0.01342 & (+)\mbox{ }0.001454 \\
          &       &       & SSO-C & ($\approx$)\mbox{ }0.012665 & (+)\mbox{ }0.012765 & (+)\mbox{ }0.012868 & (+)\mbox{ }9.29E-05 \\
          &       &       & SbPPA & 0.012665 & 0.012666 & 0.012666 & 3.39E-10 \\[0.7mm]
   17  & Welded & Not Known & PSO   & (+)\mbox{ }1.846408 & (+)\mbox{ }2.011146 & (+)\mbox{ }2.237389 & (+)\mbox{ }0.108513 \\
          &Beam Design       &       & ABC   & (+)\mbox{ }1.798173 & (+)\mbox{ }2.167358 & (+)\mbox{ }2.887044 & (+)\mbox{ }0.254266 \\
          &Problem       &       & FF    & (+)\mbox{ }1.724854 & (+)\mbox{ }2.197401 & (+)\mbox{ }2.931001 & (+)\mbox{ }0.195264 \\
          &       &       & SSO-C & ($\approx$)\mbox{ }1.724852 & (+)\mbox{ }1.746462 & (+)\mbox{ }1.799332 & (+)\mbox{ }0.02573 \\
          &       &       & SbPPA & 1.724852 & 1.724852 & 1.724852 & 4.06E-08 \\[0.7mm]
    18 & Speed & Not Known & PSO   & (+)\mbox{ }3044.453 & (+)\mbox{ }3079.262 & (+)\mbox{ }3177.515 & (+)\mbox{ }26.21731 \\
          &Reducer Design       &       & ABC   & (+)\mbox{ }2996.116 & (+)\mbox{ }2998.063 & (+)\mbox{ }3002.756 & (+)\mbox{ }6.354562 \\
          &Optimization       &       & FF    & (+)\mbox{ }2996.947 & (+)\mbox{ }3000.005 & (+)\mbox{ }3005.836 & (+)\mbox{ }8.356535 \\
          &       &       & SSO-C & ($\approx$)\mbox{ }2996.113 & ($\approx$)\mbox{ }2996.113 & ($\approx$)\mbox{ }2996.113 & (+)\mbox{ }1.34E-12 \\
          &       &       & SbPPA & 2996.114 & 2996.114 & 2996.114 & 0 \\
    \bottomrule
 \end{tabular}%
    \end{minipage}
  \label{tab:CF01}%
  \end{table}
}
\begin{figure}
\hspace{-21mm}
\includegraphics[width=19cm]{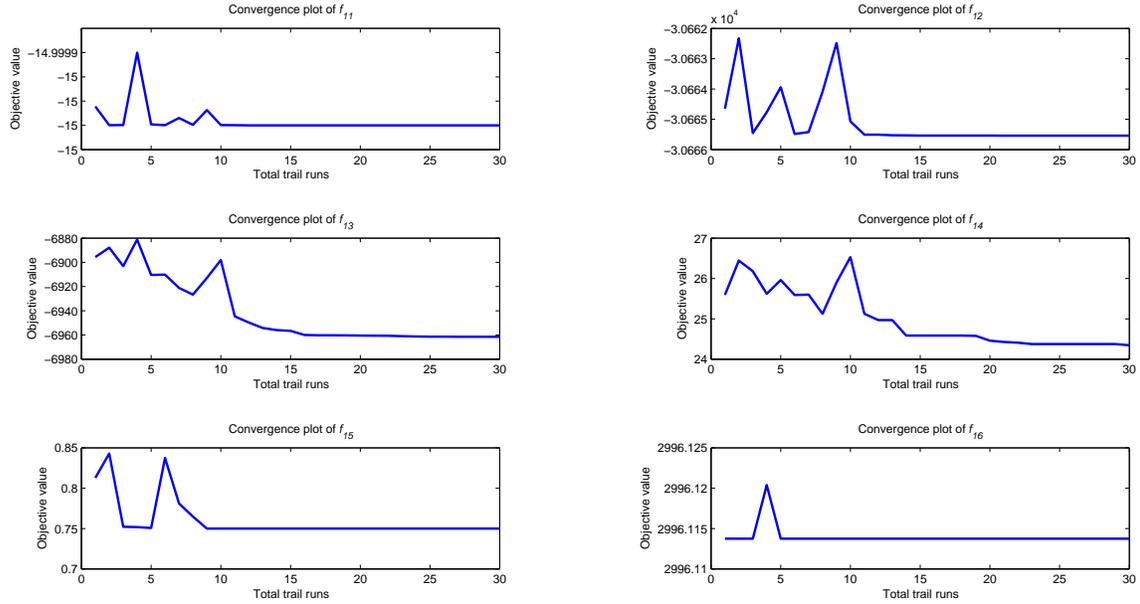}
\caption{Performance of SbPPA on constrained global optimization problems. The problems solved in this table are standard constrained optimization problems}
\end{figure}
\newpage
\section{Conclusion}
A new algorithm mimicking the  seed-based plant propagation (SbPPA) is designed and implemented for both unconstrained and constrained  optimization problems. The performance of SbPPA is compared  with a number of well established algorithms. The results are compiled in terms of best, mean, worst and standard deviation. SbPPA is very easy to implement as it needs less arbitrary parameter settings. An alternative strategy is adopted to update our current population. The effects on convergence are shown through convergence plots, Figures (4-5), of some of the solved problems. Note that the success rate of SbPPA depends on the quality of the initial population. SbPPA is being tested on discrete real world problems.
\section{Acknowledgments}
This work is supported by Abdul Wali Khan University, Mardan, Pakistan, Grant No. F.16-5/ P\& D/ AWKUM /238.

\newpage
\setcounter{section}{0}
\section{Appendix}
\setcounter{section}{0}
\section{Set of Unconstrained Global Optimization Problems}
{\fontsize{12}{12}
\selectfont
\begin{table}[htp]
 \hspace{-3mm}
  \footnotesize\setlength{\tabcolsep}{2.5pt}
  \begin{minipage}{16cm}
  \caption{Unconstrained Global Optimization Problems Used In Our Experiments.}
   \begin{tabular}{l@{\hspace{8pt}} *{7}{l}}
    \toprule
    Fun   &Ftn. Name& D     & C     & Range & Min   & Formulation \\[1mm]
   \midrule
    $f_1$& Colville & 4 & UN    & [-10  10] & 0     &$f(x)=100(x_1^2-x_2)+(x_1-1)^2+(x_3-1)^2+90(x_3^2-x_4)^2+10.1((x_2-1)^2$\\&&&&&&\hspace{8mm}$+(x_4-1)^2)+19.8(x_2-1)(x_4-1)$  \\[1mm]
    $f_2$& Matyas& 2     & UN    & [-10  10] & 0     &$f(x)=0.26(x_1^2+x_2^2)-0.48x_1x_2 $ \\[1mm]
    $f_3$&Schaffer & 2     & MN    & [-100  100] & 0     & $f(x)=0.5+\frac{\sin^{2}(\sqrt{\sum_{i=1}^{n}x^{2}_{i}})-0.5}{(1+0.001(\sum_{i=1}^{n}x^{2}_{i}))^{2}}$   \\[2mm]
    $f_4$&Six Hump Camel Back & 2     & MN    & [-5  5] & -1.03163 & $f(x)=4x_1^2-2.1x_1^4+\frac{1}{3}x_1^6+x_1x_2-4x_2^2+4x_2^4$ \\[1mm]
    $f_5$& Trid6& 6     & UN    & [-36  36] & -50   & $f(x)=\sum_{i=1}^{6} (x_i-1)^2-\sum_{i=2}^{6}x_ix_{i-1}$  \\[1mm]
    $f_6$&Trid10 & 10    & UN    & [-100  100] & -210  &$f(x)=\sum_{i=1}^{10} (x_i-1)^2-\sum_{i=2}^{10}x_ix_{i-1}$ \\[1mm]
    $f_7$&Sphere & 30    & US    & [-100  100] & 0     &$f(x)=\sum_{i=1}^{n}x^{2}_{i}$   \\[1mm]
    $f_8$&SumSquares & 30    & US    & [-10  10] & 0     & $f(x)=\sum_{i=1}^{n}ix^{2}_{i} $  \\[1mm]
    $f_9$& Griewank& 30    & MN    & [-600  600] & 0     &$f(x)=\frac{1}{4000}\sum_{i=1}^{n}x^{2}_{i}-\prod_{i=1}^{n}\cos(\frac{x_{i}}{\sqrt{i}})+1$  \\[1mm]
    $f_{10}$&Ackley & 30    & MN    & [-32  32] & 0     &$f(x)=-20\exp(-0.2\sqrt{\frac{1}{n}\sum_{i=1}^{n}x^{2}_{i}})-\exp(\frac{1}{n}\sum_{i=1}^{n}\cos(2\pi x_{i}))+20+e$  \\[1mm]
    \bottomrule
  \end{tabular}%
    \end{minipage}
  \label{tab:CF01}%
  \end{table}
}

\section{Set of Constrained Global Optimization Problems Used in Our Experiments}

\subsection{CP1}
\begin{table}[htp]
 \vspace{-4mm}
 \hspace{-2mm}
  \footnotesize\setlength{\tabcolsep}{2.5pt}
  \begin{minipage}{16cm}
   \begin{tabular}{l@{\hspace{8pt}} *{3}{l}}
                       \hspace{20mm}    &    & Min \mbox{ }\hspace{13mm}$f(x)=5\sum_{d=1}^{4}x_d-5\sum_{d=1}^{4}x_d^2-\sum_{d=5}^{13}x_d$ \hspace{15mm}\\[2mm]
                       &    &subject to \hspace{6mm}$g_1(x)=2x_1+2x_2+x_{10}+x_{11}-10\leq0$\\[2mm]
                        &   &\hspace{18.5mm}$g_2(x)=2x_1+2x_3+x_{10}+x_{12}-10\leq0$\\[2mm]
                       &    &\hspace{18.5mm}$g_3(x)=2x_2+2x_3+x_{11}+x_{12}-10\leq0$\\[2mm]
                       &    &\hspace{18.5mm}$g_4(x)=-8x_1+x_{10}\leq0$\\[2mm]
                        &   &\hspace{18.5mm}$g_5(x)=-8x_2+x_{11}\leq0$\\[2mm]
                       &    &\hspace{18.5mm}$g_6(x)=-8x_3+x_{12}\leq0$\\[2mm]
                        &   &\hspace{18.5mm}$g_7(x)=-2x_4-x_5+x_{10}\leq0$\\[2mm]
                       &    &\hspace{18.5mm}$g_8(x)=-2x_6-x_7+x_{11}\leq0$\\[2mm]
                        &    &\hspace{18.5mm}$g_9(x)=-2x_8-x_9+x_{12}\leq0$,\\[2mm]
\end{tabular}%
    \end{minipage}
  \label{tab:CF01}%
  \end{table}
\noindent  where bounds are $0 \leq x_i \leq 1\mbox{ }(i = 1, .. ., 9, 13),\mbox{ } 0 \leq x_i \leq 100\mbox{ }(i = 10, 11,
12)$. The global optimum is at $x^* = (1, 1, 1, 1, 1, 1, 1, 1, 1, , 3, 3, 3, 1),
f(x^*)= -15$.
  \newpage
  \subsection{CP2}
\begin{table}[htp]
  \vspace{-5mm}
 \hspace{-2mm}
  \footnotesize\setlength{\tabcolsep}{2.5pt}
  \begin{minipage}{16cm}
   \begin{tabular}{l@{\hspace{8pt}} *{3}{l}}
                    \hspace{20mm}    &   & Min \mbox{ }\hspace{13mm}$f(x)=5.3578547x_2 +0.8356891x_1x_5+ 37.293239x_1 - 40792.141$\\[2mm]
                         &   &subject to \hspace{6mm}$g_1(x)=85.334407+0.0056858x_2x_5+0.0006262x_1x_4-0.0022053x_3x_5-92\leq0$\\[2mm]
                        &    &\hspace{18.5mm}$g_2(x)=-85.334407-0.0056858x_2x_5-0.0006262x_1x_4+0.0022053x_3x_5\leq0$\\[2mm]
                        &    &\hspace{18.5mm}$g_3(x)=80.51249+0.0071317x2x5+0.0029955x1x2-0.0021813x_2-110\leq0$\\[2mm]
                        &    &\hspace{18.5mm}$g_4(x)=-80.51249-0.0071317x_2x_5+0.0029955x_1x_2-0.0021813x_2+90\leq0$\\[2mm]
                        &    &\hspace{18.5mm}$g_5(x)=9.300961-0.0047026x_3x_5-0.0012547x_1x_3-0.0019085x_3x_4-25\leq0$\\[2mm]
                       &     &\hspace{18.5mm}$g_6(x)=-9.300961-0.0047026x_3x_5-0.0012547x_1x_3-0.0019085x_3x_4+20\leq0$,\\[2mm]
\end{tabular}%
    \end{minipage}
  \label{tab:CF01}%
  \end{table}

\noindent where    $78 \leq x_1\leq 102$,    $33 \leq x_2 \leq 45$,    $27 \leq x_i\leq 45$ $(i = 3,    4,    5)$.
The  optimum  solution  is  $x^* = (78,  33,  29.995256025682,  45, 36.775812905788)$, where $f(x^*)= - 30665.539$. Constraints $g_1$ and
$g_6$ are active.
\subsection{CP3}
\begin{table}[htp]
 \vspace{-5mm}
 \hspace{-2mm}
  \footnotesize\setlength{\tabcolsep}{2.5pt}
  \begin{minipage}{16cm}
   \begin{tabular}{l@{\hspace{8pt}} *{3}{l}}
                    \hspace{20mm}    &   & Min \mbox{ }\hspace{13mm}$f(x)=(x_1-10)^3+(x_2-20)^3$\\[2mm]
                         &   &subject to \hspace{6mm}$g_1(x)=-(x_1-5)^2-(x_2 -5)^2+100 \leq0$\\[2mm]
                        &    &\hspace{18.5mm}$g_2(x)=(x_1-6)^2 +(x^2-5)^2-82.81 \leq0$,\\[2mm]

\end{tabular}%
    \end{minipage}
  \label{tab:CF01}%
  \end{table}
  \noindent where $13 \leq x_1 \leq 100 $ and  $0 \leq x_2 \leq 100$.  The  optimum  solution  is $x^* = (14.095, 0.84296)$ where $f(x^*)= -6961.81388$. Both constraints are active.

\subsection{CP4}
\begin{table}[htp]
  \vspace{-5mm}
 \hspace{-2mm}
  \footnotesize\setlength{\tabcolsep}{2.5pt}
  \begin{minipage}{16cm}
   \begin{tabular}{l@{\hspace{8pt}} *{3}{l}}
                    \hspace{20mm}    &   & Min \mbox{ }\hspace{13mm}$f(x)=x_1^2+x_2^2+x_1x_2-14x_1-16x_2+(x_3-10)^2+ 4(x_4 - 5)^2 + (x_5 - 3)^2 + 2(x_6 -1)^2 $\\[2mm]
                         &    &\hspace{28.5mm}$+ 5x_7^2 + 7(x_8 - 11)^2+ 2(x_9 - 10)^2 + (x_{10}-7)^2 + 45$\\[2mm]
                         &   &subject to \hspace{6mm}$g_1(x)=-105+4x_1+5x_2-3x_7+9x_8 \leq0$\\[2mm]
                        &    &\hspace{18.5mm}$g_2(x)=10x_1 - 8x_2 - 17x_7 + 2x_8 \leq0$\\[2mm]
                        &    &\hspace{18.5mm}$g_3(x)=-8x_1 + 2x_2 + 5x_9 - 2x_{10} - 12 \leq0$\\[2mm]
                        &    &\hspace{18.5mm}$g_4(x)=3(x_1 - 2)^2 + 4(x_2 - 3)^2 + 2x_3^2 - 7x_4 - 120 \leq0$\\[2mm]
                        &    &\hspace{18.5mm}$g_5(x)=5x_1^2 + 8x_2 + (x_3 - 6)^2 - 2x^4 - 40 \leq0$\\[2mm]
                       &     &\hspace{18.5mm}$g_6(x)=x_1^2 + 2(x_2 - 2)^2 - 2x_1x_2 + 14x_5 - 6x_6 \leq0$\\[2mm]
                       &     &\hspace{18.5mm}$g_7(x)=0.5(x_1 - 8)^2 + 2(x_2 - 4)^2 + 3x_5^2 - x_6 - 30  \leq0$\\[2mm]
                       &     &\hspace{18.5mm}$g_8(x)=-3x_1 + 6x_2 + 12(x_9 - 8)^2 - 7x_{10} \leq0$,\\[2mm]
\end{tabular}%
    \end{minipage}
  \label{tab:CF01}%
  \end{table}
\noindent where   $-10 \leq x_i \leq 10$ $(i = 1,   .. .,   10)$.   The   global   optimum   is\\
$x^* = (2.171996,	2.363683,	8.773926,	5.095984,	0.9906548,
1.430574,   1.321644,   9.828726,   8.280092,   8.375927)$, where $f(x^*) = 24.3062091$. Constraints $g_1,\mbox{ }g_2,\mbox{ } g_3,\mbox{ } g_4,\mbox{ } g_5$ and $g_6$ are active.

\newpage
\subsection{CP5}
\begin{table}[htp]
 \vspace{-5mm}
 \hspace{-2mm}
  \footnotesize\setlength{\tabcolsep}{2.5pt}
  \begin{minipage}{16cm}
   \begin{tabular}{l@{\hspace{8pt}} *{3}{l}}
                    \hspace{20mm}    &   & Min \mbox{ }\hspace{13mm}$f(x)=x_1^2 + (x_2 - 1)^2$\\[2mm]
                         &   &subject to \hspace{6mm}$g_1(x)=x_2 - x_1^2=0$,\\[2mm]

\end{tabular}%
    \end{minipage}
  \label{tab:CF01}%
  \end{table}

\noindent where  $−1 \leq x_1 \leq 1$,  $−1 \leq x_2 \leq 1$.  The  optimum  solution  is  $x^*=
(\pm1/\sqrt{(2)}, 1/2)$, \\where $f(x^*) = 0.7499$.

\subsection{Welded Beam Design Optimisation}
 The welded beam design is a standard test problem for constrained design optimisation \cite{cagnina2008solving,yang2010engineering}. There are four design variables: the width $w$ and length $L$ of the welded area, the depth $d$ and thickness $h$ of the main beam. The objective is to minimise the overall fabrication cost, under the appropriate constraints of shear stress $\tau$, bending stress $\sigma$, buckling load $P$ and maximum end deflection $\delta$.
The optimization model is summarized as follows, where $x^T=(w,L,d,h).$
\begin{minipage}{12cm}
\begin{equation}
\hspace{-7mm} Minimise\hspace{5mm}f(x) = 1.10471w^2L + 0.04811dh(14.0 + L),
\end{equation}
\hspace{15mm}subject to
\begin{eqnarray}
\begin{aligned}
  &g_{1}(x) = w-h \leq 0,  \\
  &g_{2}(x) = \delta(x)-0.25 \leq 0,  \\
  &g_{3}(x) = \tau(x)-13,600 \leq 0, \\
  &g_{4}(x) =  \sigma(x)-30,000\leq 0, \\
  &g_{5}(x) = 1.10471w^2 + 0.04811dh(14.0 + L)-5.0 \leq 0,  \\
  &g_{6}(x) =  0.125-w\leq 0,  \\
  &g_{7}(x) = 6000-P(x) \leq 0,\\
  \end{aligned}
\end{eqnarray}
\end{minipage}\\
where
\begin{equation}
\begin{aligned}
&\sigma(x) =\frac{504,000}{hd^2},\\[2mm]
&D = \frac{1}{2}\sqrt{L^2+(w+d)^2},\\[2mm]
&\delta =\frac{65,856}{30,000hd^3},\\[2mm]
&\alpha = \frac{6000}{\sqrt{2}wL},\\[2mm]
&P =0.61423 \times 10^6 \frac{dh^3}{6}\left(1-\frac{\sqrt[d]{\frac{30}{48}}}{28}\right).\\
\end{aligned}
\begin{aligned}
&Q =6000\left(14+\frac{L}{2}\right),\\[2mm]
&J = \sqrt{2}w L \left(\frac{L^2}{6}+\frac{(w+d)^2}{2}\right),\\[2mm]
&\beta = \frac{QD}{J},\\[2mm]
\indent \indent&\tau(x) =\sqrt{\alpha^2+\frac{\alpha\beta L}{D}+\beta^2}.\\
\end{aligned}
\end{equation}
\subsection{Speed Reducer Design Optimization}
The problem of designing a speed reducer \cite{golinski1974adaptive} is a standard test problem. It consists of the design variables as: face width $x_1$, module of teeth $x_2$, number of teeth on pinion $x_3$, length of the first shaft between bearings $x_4$, length of the second shaft between bearings $x_5$, diameter of the first shaft $x_6$, and diameter of the first shaft $x_7$ (all variables continuous except $x_3$ that is integer). The weight of the speed reducer is to be minimized subject to constraints on bending stress of the gear teeth, surface stress, transverse deflections of the shafts and stresses in the shaft, \cite{cagnina2008solving}. The mathematical formulation of the problem, where $x^T=(x_1,x_2,x_3,x_4,x_5,x_6,x_7)$, is as follows.
\vspace{-5mm}
\begin{multline}
\begin{aligned}
 Minimise \indent f(x) = &0.7854x_1x_2^2(3.3333x_3^2 + 14.9334x_3 − 43.0934)\\
 &-1.508x_1(x_6^2 + x_7^3 ) + 7.4777(x_6^3 + x_7^3 )+ 0.7854(x_4x_6^2 + x_5x_7^2 ),
\end{aligned}
\end{multline}
\vspace{-5mm}
\hspace{5mm}subject to
\vspace{-.5mm}
\begin{eqnarray}
\begin{aligned}
  &g_{1}(x) = \frac{27}{x_1x_2^2x_3}-1\leq 0,  \\
  &g_{2}(x) = \frac{397.5}{x_1x_2^2x_3^2}-1 \leq 0,  \\
  &g_{3}(x) = \frac{1.93x_4^3}{x_2x_3x_6^4}-1 \leq 0, \\
  &g_{4}(x) =  \frac{1.93x_5^3}{x_2x_3x_7^4}-1\leq 0, \\
  &g_{5}(x) = \frac{1.0}{110x_6^3}\sqrt{\left(\frac{745.0x_4}{x_2x_3}\right)^2+16.9\times10^6}-1 \leq 0,  \\
  &g_{6}(x) = \frac{1.0}{85x_7^3}\sqrt{\left(\frac{745.0x_5}{x_2x_3}\right)^2+157.5\times10^6}-1 \leq 0,  \\
  &g_{7}(x) = \frac{x_2x_3}{40}-1 \leq 0,\\
  &g_{8}(x) =  \frac{5x_2}{x_1}-1\leq 0, \\
  &g_{9}(x) = \frac{x_1}{12x_2}-1 \leq 0,  \\
  &g_{10}(x) =  \frac{1.5x_6+1.9}{x_4}-1\leq 0,  \\
  &g_{11}(x) = \frac{1.1x_7+1.9}{x_5}-1 \leq 0.\\
  \end{aligned}
\end{eqnarray}
The simple limits on the design variables are \\ $2.6\leq x_1\leq3.6 $, $0.7 \leq x_2 \leq0.8$, \\ $17\leq x_3\leq28$, $7.3\leq x_4\leq8.3 $, $7.8 \leq x_5 \leq8.3$, \\ $2.9\leq x_6\leq3.9$ and $5.0\leq x_7\leq5.5$.
\subsection{Spring Design Optimisation}
The main objective of this problem \cite{arora2004introduction,belegundu1985study} is to minimize the weight of a tension/compression string, subject to constraints of minimum deflection, shear stress, surge frequency, and limits on outside diameter and on design variables. There are three design variables: the wire diameter $x_1$, the mean coil diameter $x_2$, and the number of active coils $x_3$, \cite{cagnina2008solving}. The mathematical formulation of this problem, where $x^T=(x_1,x_2,x_3)$, is as follows.
\begin{eqnarray}
\begin{aligned}
 \hspace{-25.5mm}Minimize \mbox{ }f(x)= (x_3+2)x_2x_1^2,
\end{aligned}
\end{eqnarray}

\hspace{35.5mm}subject to
\vspace{-5mm}
\begin{eqnarray}
\begin{aligned}
  &g_{1}(x) = 1-\frac{x_2^3x_3}{7,178x_1^4}\leq 0,  \\
  &g_{2}(x) = \frac{4x_2^2-x_1x_2}{12,566(x_2x_1^3)-x_1^4}+\frac{1}{5,108x_1^2}-1 \leq 0,  \\
  &g_{3}(x) = 1-\frac{140.45x_1}{x_2^2x_3} \leq 0, \\
  &g_{4}(x) = \frac{x_2+x_1}{1.5}-1\leq 0. \\
\end{aligned}
\end{eqnarray}
The simple limits on the design variables are $0.05\leq x_1\leq 2.0 $, $0.25 \leq x_2 \leq1.3$ \\and $2.0\leq x_3\leq15.0.$

\newpage

\bibliography{Sulaimanbib}
\bibliographystyle{plain}
\end{document}